\documentclass[english,draft]{article}
\usepackage[english]{babel}
\usepackage{amsmath}
\usepackage{amsthm}
\usepackage{amssymb}
\usepackage{color}

\newtheorem{theorem}{Theorem}[]
\newtheorem{lemma}{Lemma}[]

\newtheorem{problem}{Problem}[]
\newtheorem*{remark}{Remark}

\mathversion{bold}
\begin{document}
\title{The concrete theory of numbers : Problem of simplicity of Fermat number-twins}
\author{Boris\, V. Tarasov\thanks{Tarasov Boris V.
The concrete theory of numbers: Problem of simplicity of Fermat number-twins.
MSC 11A51.
\textcolor[rgb]{1.00,0.00,0.50}{\textcopyright 2007 Tarasov B.\,V.,}
\textcolor[rgb]{1.00,0.00,0.50}{independent researcher of Unknown.}
}\\
--------------------------------------\\
}
\maketitle
\begin{abstract}
The problem of simplicity of Fermat number-twins \\$f_{n}^{\pm}=2^{2^n}\pm3$ is studied.
The question for what $n$ numbers $f_{n}^{\pm}$ are composite is investigated.
The factor-identities for numbers of a kind $x^2 \pm k $ are found.
\end{abstract}

{\centering\section{Introduction}\par}

In present work we consider Fermat numbers
\begin{equation}\label{e:NumF}f_{n}=2^{2^n}+1,\end{equation}
where $n\geq0$ is integer.\par
Fermat number-twins we will define as :
\begin{equation}\label{e:NumTwinF}
 \begin{aligned}
  f_{n}^{-} &= f_{n}-4 = 2^{2^n}-3,\\
  f_{n}^{+} &= f_{n}+2 = 2^{2^n}+3,\\
 \end{aligned}
\end{equation}
where $n\geq0$  is integer( for $f_{n}^{-}$ is considered $n\geq1$ ).\par
Pierre de Fermat in 1640 in the letter to Mersenne \cite{Fermat} suggested a hypothesis,
that all numbers (\ref{e:NumF}) are
prime, which has been denied by Leonhard Euler \cite{Euler} in 1732, Euler has found decomposition
\begin{displaymath} f_{5}=641\cdot6700417. \end{displaymath}\par
Eisenstein (1844) assumed, that there exists infinite number of prime Fermat numbers.
There was a well-known problem of Fermat prime \cite[Eisenstein]{Arnold,Weisstein}.\par
\begin{problem}[\textbf{\textcolor[rgb]{0.00,0.00,1.00}{Fermat prime}}]~Whether there are infinitely many\\
Fermat prime ?\end{problem}
There was also other not less known problem :
\begin{problem}[\textbf{\textcolor[rgb]{0.00,0.00,1.00}{Fermat composite numbers}}]
~Whether there are infinitely many Fermat composite numbers ?\end{problem}
Problems 1 and 2 are actual even for today. Therefore the researches, throwing light on numbers,
located closely to Fermat numbers are extremely useful. Probably, the similar information
will reveal a secret of the Fermat numbers in the future.\par
Five Fermat prime numbers \cite[Fermat prime]{Arnold,Graham,Weisstein,SearchFermatNumberDivisors} are known
$f_{0}=3$, $f_{1}=5$, $f_{2}=17$, $f_{3}=257$, $f_{4}=65537$.
The author knows only five Fermat prime number-twins( see \cite{Sloane} ) : \\
$f_{2}^{-}=13$, $f_{0}^{+}=5$, $f_{1}^{+}=7$, $f_{2}^{+}=19$, $f_{4}^{+}=65539$.\par

{\centering\section{Prime divisors of Fermat number-twins}\par}

Let's bring the simple statements concerning Fermat number-twins.
\begin{theorem}[] The following congruences are fair :\par
(1) If $n$ is an even number, then
\begin{equation}\label{e:M1}f_{n}^{-}\equiv0(mod\,13).\end{equation}\par
(2) If $n=4k+3$, where $k\geq0$, then
\begin{equation}\label{e:M2}f_{n}^{-}\equiv0(mod\,11).\end{equation}
\end{theorem}
\begin{proof}
1) Let $n=2$, then $f_{2}^{-}=2^{4}-3=13$. Let's assume, that congruence (\ref{e:M1}) is proved
for all even numbers $n$, where $n\leq k$, $k$ is even number. Let's consider expression\\
$f_{k+2}^{-} - f_{k}^{-} = 2^{2^{k+2}} - 2^{2^{k}} = (2^{2^{k}})^4 - 2^{2^{k}} = 2^{2^{k}}\{(2^{2^{k}})^3 - 1\} =$\\
$ = 2^{2^{k}}(2^{2^{k}} - 1)\{(2^{2^{k}})^2 + 2^{2^{k}} + 1\}$. As $2^{2^{k}}\equiv3(mod\,13)$, that\\
$(2^{2^{k}})^2 + 2^{2^{k}} + 1 \equiv 3^2 + 3 + 1 = 13 \equiv0(mod\,13) $. We have proved, that \\
$f_{k+2}^{-}\equiv0(mod\,13)$.\par
2) If $k = 0$, then $n=3$, $f_{3}^{-}=2^{2^3}-3=2^{8}-3=11\cdot23$. \\ Let's consider expression
$f_{4k+3}^{-} - f_{3}^{-} = 2^{8\cdot2^{4k}} - 2^{8} = $ \\
$ = 2^{8}\{(2^{8})^{2^{4k}-1}-1\}$. Further, as $2^{5}\equiv-1(mod\,11)$, \\$2^{15}\equiv-1(mod\,11)$ and
$2^{4k}-1\equiv0(mod\,15)$, that\\ $(2^{8})^{2^{4k}-1}-1\equiv0(mod\,11)$.
\end{proof}
\begin{theorem}[] If $n$ is an odd number, then
\begin{equation}\label{e:M3}f_{n}^{+}\equiv0(mod\,7).\end{equation}
\end{theorem}
\begin{proof}
$f_{1}^{+}=7$. Let's assume, that congruence (\ref{e:M3}) is proved for all odd numbers
$n$, where $n\leq k$, $k$ is odd number. Let's consider expression\\
$f_{k+2}^{+} - f_{k}^{+} = 2^{2^{k+2}} - 2^{2^{k}} = (2^{2^{k}})^4 - 2^{2^{k}} = 2^{2^{k}}\{(2^{2^{k}})^3 - 1\} =$\\
$ = 2^{2^{k}}(2^{2^{k}} - 1)\{(2^{2^{k}})^2 + 2^{2^{k}} + 1\}$. As $2^{2^{k}}\equiv-3(mod\,7)$, that\\
$(2^{2^{k}})^2 + 2^{2^{k}} + 1 \equiv 3^2 - 3 + 1 = 7 \equiv0(mod\,7) $. We have proved, that \\
$f_{k+2}^{-}\equiv0(mod\,7)$.
\end{proof}
\begin{lemma}[] The following congruences are fair :\par
\begin{equation}\label{e:M4}
 \begin{aligned}
     2^{18k+2}+3 &\equiv0(mod\,7),\\
     2^{18k+4}+3 &\equiv0(mod\,19),\\
     2^{18k+8}+3 &\equiv0(mod\,7),\\
     2^{18k+14}+3 &\equiv0(mod\,7),\\
 \end{aligned}
\end{equation}
where $k\geq0$ is integer.
\end{lemma}
\begin{proof}
Validity of congruences (\ref{e:M4}) follows obviously from the equality\\
$2^{18}-1=3^{3}\cdot7\cdot19\cdot73$.
\end{proof}
\begin{lemma}[] The following statements are fair, where $t\geq0$ is integer :\par
(1) If $2^{t}-1\equiv0(mod\,9)$, then $t=6k$.\par
(2) If $2^{t}+1\equiv0(mod\,9)$, then $t=6k+3$.
\end{lemma}
\begin{proof}
1) Let $2^{t}-1\equiv0(mod\,9)$ and $t=6k+r$, where\\
$0\leq r <6$. As $2^{6}-1\equiv0(mod\,9)$, that $2^{t}-1\equiv2^{r}-1\equiv0(mod\,9)$.
As $2^{0}-1\equiv0(mod\,9)$, but $2^{1}-1=1\not\equiv0(mod\,9)$, $2^{2}-1=3\not\equiv0(mod\,9)$,
$2^{3}-1=7\not\equiv0(mod\,9)$, $2^{4}-1=15\not\equiv0(mod\,9)$, $2^{5}-1=31\not\equiv0(mod\,9)$,
that $r=0$. We have proved, that $t=6k$.\\
2) Let $2^{t}+1\equiv0(mod\,9)$ and $t=6k+r$, where $0\leq r <6$.\\
Then $2^{t}+1\equiv2^{r}+1\equiv0(mod\,9)$. As $2^{0}+1=2\not\equiv0(mod\,9)$,\\
$2^{1}+1=3\not\equiv0(mod\,9)$, $2^{2}+1=5\not\equiv0(mod\,9)$, $2^{4}+1=17\not\equiv0(mod\,9)$,
$2^{5}+1=33\not\equiv0(mod\,9)$, but $2^{3}+1=9\equiv0(mod\,9)$, that $r=3$. We have proved, that $t=6k+3$.
\end{proof}
\begin{theorem}[] The following statements are fair :\par
(1) Prime numbers $f_{n}^{-}$, $n>2$ are possible only for $n=4k+1$, where $k\geq1$. \par
(2) Prime numbers $f_{n}^{+}$, $n>4$ are possible only for $n=6k$ or $n=6k+4$, where $k\geq1$.
\end{theorem}
\begin{proof}
The statement (1) is the corollary of the theorem 1. Let's prove \\the statement (2). Let's assume,
that number $f_{n}^{+}$ is prime.
Then it follows from a lemma 1, that if $2^{n}=18m+r$, where $0\leq r <18$, $r$ is an even number, then $r=10$ or $r=16$.
If $2^{n}=18m+10$, then $2^{n}-1\equiv0(mod\,9)$. Then it follows from a lemma 2, that $n=6k$. If $2^{n}=18m+16$, then
$2^{n}+2\equiv0(mod\,9)$ or $2^{n-1}+1\equiv0(mod\,9)$. Then it follows from a lemma 2, that $n-1=6k+3$ or $n=6k+4$.
\end{proof}\par
It is checked up, that for $n>2$, $n=4k+1$, where $0< k \leq4$, i.e. for $n=5, 9, 13, 17$  numbers $f_{n}^{-}$ are
composite.
It is checked up, that for $n>4$, $n=6k$ and $n>4$, $n=6k+4$, where $0< k \leq2$, i.e. for $n=6, 10, 12, 16$
numbers $f_{n}^{+}$ are composite.
\begin{remark}[MMonline] \textcolor[rgb]{0.98,0.00,0.00}{M. Alekseyev} \cite{Tarasov} at the forum MMonline\\
"Mathematics"\ furnishs without the proof the statement :\\
\textcolor[rgb]{0.00,0.00,1.00}{If} $\textcolor[rgb]{0.00,0.00,1.00}{n=12k+10}$\textcolor[rgb]{0.00,0.00,1.00}{, then}
$\textcolor[rgb]{0.00,0.00,1.00}{f_{n}^{+}\equiv0(mod\,79)}$\textcolor[rgb]{0.00,0.00,1.00}{.}
\end{remark}
Thus, $f_{n}^{+}$ are composite and for $n=6k+4$, where $k\geq1$ is an odd number.
\begin{theorem}[] Numbers $f_{n}^{-}$ and $f_{n}^{+}$ are \textbf{\textcolor[rgb]{0.00,0.00,1.00}{simultaneously composite}}\\
for all $n$ of kind :\par
\begin{equation}\label{e:M5}
 \begin{aligned}
 n &= 12k + 2,     &(1^{\circ})\\
 n &= 12k + 3,     &(2^{\circ})\\
 n &= 12k + 7,     &(3^{\circ})\\
 n &= 12k + 8,     &(4^{\circ})\\
 n &= 12k + 11,     &(5^{\circ})\\
 \end{aligned}
\end{equation}
where $k\geq1$ is integer.
\end{theorem}
\begin{proof}
Validity of the theorem 4 follows obviously from statements of the theorem 3.
\end{proof}

{\centering\section{The factor-identities for numbers of a kind $X^2\pm k$}\par}
For the beginning, let's bring two known identities.

\subsection{The first factor-identity}
\begin{equation}\label{e:Id1}1 + (L^{2} + L + 1)^{2} =
(L^{2} + 1)\cdot(L^{2} + 2L + 2),\end{equation}
where $L\geq1$ is integer.

\subsection{The second factor-identity}
\begin{equation}\label{e:Id2}1 + (2L^{2})^{2} =
(2L^{2} - 2L + 1)\cdot(2L^{2} + 2L + 1),\end{equation}
where $L > 1$ is integer.
\pagebreak
\subsection{The new factor-identity}
If\\
\begin{equation}\label{e:Idnew}
 \begin{aligned}
  A &= a + amnk + n^{2}/4\cdot(am^{2} + n^{2})k^{2},\\
  B &= 1 - mnk  + m^{2}/4\cdot(am^{2} + n^{2})k^{2},\\
  X &= 1/2(am^{2} - n^{2})k + mn/4\cdot(am^{2} + n^{2})k^{2},\\
 \end{aligned}
\end{equation}
where $m$,$n$,$k$,$a$ are integer, then equality takes place
\begin{equation} a + X^{2} = A\cdot B.\end{equation}
\par The author does not know the factor-identity of type $1 + X^2 = AB$, left part of which
contains an infinite subset of Fermat numbers. Existence of such factor-identity positively solves the
\textbf{\textcolor[rgb]{0.00,0.00,1.00}{problem 2}}, that there exists the infinite number of composite Fermat numbers !

{\centering\section{The conclusion}\par}
If number of prime Fermat $f_{n}$ is finite, that, starting from a number $k_{0} > 1$, for all numbers
 $n$ from (\ref{e:M5}), where $k > k_{0}$, we receive an infinite sequence of compound numbers with 9 numbers in a part,
 namely
\begin{equation}\label{e:M6}f_{n} - 5, f_{n} - 4, f_{n} - 3, f_{n} - 2, f_{n} - 1,
f_{n}, f_{n} + 1, f_{n} + 2, f_{n} + 3.\end{equation}\par
Such gaps from composite numbers are too regular. The author supports
\textcolor[rgb]{1.00,0.00,0.50}{the assumption}
of the Eisenstein(1844),
that \textcolor[rgb]{1.00,0.00,0.50}{there exists the infinite number of prime Fermat numbers !}\par
The author offers a problem as a unresolved task :
\begin{problem}[\textbf{\textcolor[rgb]{0.00,0.00,1.00}{Prime Fermat number-twins}}]
~Whether there are prime\\ Fermat number-twins $f_{n}^{\pm}$ for
$n > 4$ ?
\end{problem}\par
As it was already noted, numbers $f_{2}^{-}=13$, $f_{0}^{+}=5$, $f_{1}^{+}=7$, $f_{2}^{+}=19$, $f_{4}^{+}=65539$ are prime.
It is checked up by the author, that for $n \leq 17$ there is no other prime Fermat number-twins.

\pagebreak

\par

\textcolor[rgb]{0.00,0.00,1.00}{---------------------------------------------------------------------}\par
\textcolor[rgb]{0.00,0.00,1.00}{Institute of Thermophysics, Siberian Branch of RAS }\par
\textcolor[rgb]{0.00,0.00,1.00}{Lavrentyev Ave., 1, Novosibirsk, 630090, Russia }\par
\textcolor[rgb]{0.00,0.00,1.00}{E-mail: tarasov@itp.nsc.ru }\par
\textcolor[rgb]{0.00,0.00,1.00}{---------------------------------------------------------------------}\par
\textcolor[rgb]{1.00,0.00,0.50}{Boris Vladimirovich Tarasov, }\par
\textcolor[rgb]{1.00,0.00,0.50}{Independent researcher of Unknown, }\par
\textcolor[rgb]{1.00,0.00,0.50}{E-mail: tarasov-b@mail.ru }\par
\textcolor[rgb]{0.00,0.00,1.00}{---------------------------------------------------------------------}\par


\begin{thebibliography}{9}
\par
\bibitem{Fermat} Pierre de Fermat.
\textit{Letter to Marin Mersenne(25 December 1640), CEuvres de Fermat, volume 2, 212-217.}
\bibitem{Euler} Euler, L.
\textit{"Observationes de theoremate quodam Fermatiano aliisque ad numeros primos spectantibus."
Acad. Sci. Petropol. 6, 103-107, ad annos 1732-33(1738). In Leonhardi Euleri Opera Omnia, Ser. I, Vol. II.
Leipzig: Teubner, pp. 1-5, 1915.}
\bibitem{Arnold} Arnold I.V.
\textit{Teoriya chisel.\,-\,M.\,:\,Uchpedgiz,\,1939.}
\bibitem{Vinogradov} Vinogradov I.\,M.
\textit{Osnovy teorii chisel.\,-\,M.\,:\,Nauka,\,1981.}
\bibitem{Graham} Ronald L.\,Graham,\,Donald E.\,Knuth,\,Oren Patashnik.\\
\textit{Concrete Mathematics :\,A Foundation for Computer Science,\,2nd edition
(Reading,\,Massachusetts:\,Addison-Wesley), 1994.}
\bibitem{Ribenboim} Ribenboim,\,P.
\textit{"Fermat Numbers"\,and\,"Numbers $k\times2^n\pm1$." \textsection2.6 and 5.7 in The New Book of Prime Number Records.
New York: Springer-Verlag, pp. 83-90 and 355-360, 1996.}
\bibitem{Weisstein} Weisstein,\,Eric W.
\textit{"Fermat Number". From MathWorld--A Wolfram Web Resource.
---http://mathworld.wolfram.com/FermatNumber.html/.\\
\textcopyright 1999---2007 Wolfram Research,\,Inc.}
\bibitem{SearchFermatNumberDivisors}
\textit{Distributed Search for Fermat Number Divisors.\\
---http://www.fermatsearch.org/.}
\bibitem{Sloane} Sloane,\,N.\,J.\,A.
\textit{Sequences A057732 and A050414 in "The On-Line Encyclopedia of Integer Sequences."}
\bibitem{Tarasov} Boris Tarasov.
\textit{Interesnaya zadachka na sravnimost' celix chisel !\\
---http://www.mmonline.ru/forum/read.php?f=1\&i=6616\&t=6616}
\end{thebibliography}
\end{document}